\documentclass{amsart}
\usepackage{amscd}

\newcommand{\w}{\omega}

\newcommand{\diam}{\mathrm{diam}\,}

\newcommand{\mesh}{\mathrm{mesh}\,}

\newcommand{\IN}{\mathbb{N}}

\newcommand{\C}{\mathcal{C}}
\newcommand{\e}{\varepsilon}
\newcommand{\Ra}{\Rightarrow}

\newcommand{\U}{\mathcal U}
\newcommand{\V}{\mathcal V}

\newcommand{\lev}{\mathrm{lev}}
\newcommand{\Deg}{\mathrm{Deg}}

\newcommand{\suc}{\mathrm{pred}}

\newcommand{\Ent}{\operatorname{Cov}}
\newcommand{\ent}{\operatorname{cov}}
\newcommand{\Cov}{\operatorname{Cov}}
\newcommand{\cov}{\operatorname{cov}}

\newcommand{\upa}{\uparrow}
\newcommand{\da}{\downarrow}

\newcommand{\Lev}{\mathrm{Lev}}

\newtheorem{theorem}{Theorem}
\newtheorem{proposition}{Proposition}
\newtheorem{lemma}{Lemma}

\theoremstyle{definition}

\newtheorem{definition}{Definition}

\title{A coarse characterization of the Baire macro-space}
\author{T.~Banakh, I.~Zarichnyi}
\address{Department of Mathematics, Ivan Franko National University of Lviv, Ukraine}
\address{Institute of Applied Problems of Mechanics and Mathematics of Ukrainian Academy of Sciences, Lviv, Ukraine}
\email{tbanakh@yahoo.com, ihor.zarichnyj@gmail.com}

\begin{document}
\begin{abstract} We prove that each coarsely homogenous separable metric space $X$ is coarsely equivalent to one of the spaces: the sigleton $1$, the Cantor macro-cube $2^{<\IN}$ or the Baire macro-space $\w^{<\IN}$. This classification is derived from coarse characterizations of the Cantor macro-cube $2^{<\IN}$ given in \cite{BZ} and of the Baire macro-space $\w^{<\IN}$ given in this paper. Namely, we prove that a separable metric space $X$ is coarsely equivalent to $\w^{<\IN}$ if any only if $X$ has asymptotic dimension zero and has unbounded geometry in the sense that for every $\delta<\infty$ there is $\e<\infty$ such that no $\e$-ball in $X$ can be covered by finitely many sets of diameter $\le \delta$.
\end{abstract}
\maketitle

This paper is devoted to the characterization of the Baire macro-space in the coarse category.
The Baire macro-space $\w^{<\IN}$ is an asymptotic counterpart of the classical Baire space $\w^\w=\prod^\w\w$ which is the Tychonoff product of countably many copies of $\w$. The Baire macro-space is defined as the countable coproduct of countably many copies of $\w$. For a non zero cardinal  $\kappa$ the  coproduct
$$\kappa^{<\IN}=\coprod\limits_{i\in\IN}\kappa=\{(x_i)_{i\in\IN}\in\kappa^\IN\colon\exists n\in \IN\;\;\forall i>n \;\; x_i=0\}$$is a metric space endowed with the ultrametric
$$d((x_i)_{i\in\IN},(y_i)_{i\in\IN})=\max(\{0\}\cup\{i\in\IN\colon x_i\neq y_i\}).$$

For $\kappa=2$ and $\kappa=\omega$ the  coproducts
$\kappa^{<\IN}$ have special names:
\begin{itemize}
\item $2^{<\IN}$ is called the {\em Cantor macro-cube};
\item $\w^{<\IN}$ is called the {\em Baire macro-space}.
\end{itemize}
In Theorem~\ref{kl-odn} we shall prove that up to the coarse equivalence these two spaces exhaust all possible types of coarsely homogeneous unbounded separable metric spaces of asymptotic dimension zero.

The coarse equivalence of metric spaces can be defined with help of multi-maps.
By a {\em multi-map} $\Phi:X\Ra Y$ between two sets $X,Y$ we
understand any subset $\Phi\subset X\times Y$. For a subset $A\subset X$ by $\Phi(A)=\{y\in Y:\exists a\in
A\mbox{ with }(a,y)\in\Phi\}$ we denote the image of $A$ under the
multi-map $\Phi$. Given a point $x\in X$ we write $\Phi(x)$
instead of $\Phi(\{x\})$.

The inverse $\Phi^{-1}:Y\Ra X$ of the multi-map $\Phi$ is the
multi-map $$\Phi^{-1}=\{(y,x)\in Y\times X: (x,y)\in\Phi\}\subset
Y\times X$$ assigning to each point $y\in Y$ the set $\Phi^{-1}(y)=\{x\in X:y\in\Phi(x)\}$. For two multi-maps $\Phi:X\Ra Y$ and $\Psi:Y\Ra Z$ we
define their composition $\Psi\circ\Phi:X\Ra Z$ as usual:
$$\Psi\circ\Phi=\{(x,z)\in X\times Z:\exists y\in Y\mbox{ such that $(x,y)\in \Phi$ and $(y,z)\in\Psi$}\}.$$

A multi-map $\Phi$ is called {\em surjective} if $\Phi(X)=Y$ and {\em bijective} if $\Phi\subset X\times Y$ coincides with the graph of a bijective (single-valued) function.

The {\em oscillation} of a multi-map $\Phi:X\Ra Y$ between metric
spaces is the function $\w_\Phi:[0,\infty)\to[0,\infty]$ assigning
to each $\delta\ge 0$ the (finite or infinite) number
$$\w_\Phi(\delta)=\sup\{\diam(\Phi(A)):A\subset
X,\;\diam(A)\le\delta\}.$$ Observe that $\w_\Phi(0)=0$ if and
only if $\Phi$ is at most single-valued in the sense that
$|\Phi(x)|\le 1$ for any $x\in X$.

A multi-map $\Phi:X\Ra Y$ between metric spaces $X$ and $Y$ is
called {\em macro-uniform} if for every $\delta<\infty$ the oscillation $\w_\Phi(\delta)$ is finite.

A multi-map $\Phi:X\Ra Y$  is called a {\em macro-uniform} {\em embedding} if
$\Phi^{-1}(Y)=X$ and both multi-maps $\Phi$ and $\Phi^{-1}$ are macro-uniform. If, in addition,
$\Phi(X)=Y$, then $\Phi$ is called  a {\em macro-uniform} {\em equivalence}.
Two metric spaces $X,Y$ are called {\em macro-uniformly} {\em equivalent\/}  if
there is a macro-uniform
equivalence $\Phi:X\Ra Y$.

Let $\e\in[0,\infty)$. By the {\em $\e$-connected component} of a point $x$ of a metric space $X$ we understand the subset $C_\e(x)$ consisting of all points $y\in X$ that can be linked with $x$ by a sequence of points $x=x_0,\dots,x_n=y$ such that $d(x_{i-1},x_i)\le\e$ for all $i\le n$. Such a sequence $x_1,\dots,x_n$ is called an {\em $\e$-chain}. It is easy to see that two $\e$-connected components $C_\e(x)$, $C_\e(y)$ either coincide or are {\em $\e$-disjoint} in the sense that $d(x',y')>\e$ for any points $x'\in C_\e(x)$, $y'\in C_\e(y)$. Thus, $\C_\e(X)=\{C_\e(x):x\in X\}$ is a disjoint cover of the metric space $X$.

In an ultrametric space $X$ the $\e$-connected components $C_\w(x)$ coincide with closed $\e$-balls $B_\e(x)=\{y\in X:d(x,y)\le\e\}$. We recall that a metric space $X$ is {\em ultrametric} if its metric $d$ satisfies the strong triangle inequality:
$$d(x,z)\le\max\{d(x,y),d(y,z)\}\mbox{ \ for all $x,y,z\in X$}.$$

A metric space $X$ has {\em asymptotic dimension zero } if  for all $\e>0$ the cover $\C_\e(X)$ has $\mesh\C_\e(X)=\sup_{x\in X}\diam C_\e(x)<\infty$. It is known that each metric space $X$ of asymptotic dimension zero
 is macro-uniformly equivalent to an ultrametric space
\cite{BDHM}.

Next, we need to introduce two cardinal characteristics $\cov_\delta^\e(X)$ and $\Cov_\delta^\e(X)$ of a metric space $X$ related to capacities of it's balls.
For a subset $A\subset X$ let $\ent_\delta(A)$ be the smallest cardinality $|\mathcal U|$ of cover $\mathcal U$ of $A$ with $\mesh(\mathcal U)\le \delta$, where $\mesh(\mathcal U)=\sup_{U\in\mathcal U}\diam U$.

For positive real numbers $\delta,\e$ consider the following two cardinals:
$$
\ent_\delta^\e(X)=\min_{x\in X} \cov_\delta (B_\e(x)) \mbox{ \ and \ }
\Ent_\delta^\e(X)=\sup_{x\in X} \ent_\delta (B_\e(x)),
$$
where $B_\e(x)$ stands for the closed $\e$-ball centered at $x$.

\begin{definition}
We say that a metric space  $X$
\begin{itemize}
\item has {\em bounded geometry, } if  there exists $\delta<\infty$ such that    $\Cov^\e_\delta(X)<\infty$ for every $\e<\infty$;
\item has {\em unbounded geometry, } if for every
$\delta<\infty$ there  exists $\e<\infty$ such that $\cov^\e_\delta(X)\ge \w$;
\item has {\em asymptotically isolated balls} if there is $\delta<\infty$ such that for every $\e<\infty$ $\cov_\delta^\e(X)=1$.
\end{itemize}
\end{definition}

Finally we recall the definition of a coarsely homogeneous metric space, introduced and studied in \cite{BCL}.

A metric space $X$ is called
\begin{itemize}
\item {\em isometrically homogeneous} if for any points $x,y\in X$
there is a bijective isometry $f:X\to X$ such that $f(x)=y$; \item
{\em coarsely homogeneous} if there is a function
$\varphi:[0,\infty)\to[0,\infty)$ such that for any points $x,y\in
X$ there is a macro-uniform equivalence $\Phi:X\Ra X$ such that
$y\in\Phi(x)$ and $\w_\Phi\le\varphi$ and
$\w_{\Phi^{-1}}\le\varphi$.
\end{itemize}

It is clear that each isometrically homogeneous metric space is coarsely homogeneous. In particular, for each cardinal $\kappa$ the space $\kappa^{<\IN}$ is isometrically and coarsely homogeneous. By \cite{BCL}, the coarse homogeneity is preserved by macro-uniform equivalences. So, each metric space that is coarsely equivalent to $2^{<\IN}$ or $\w^{<\IN}$ is coarsely homogeneous.

\begin{theorem}[Macro-classification]\label{kl-odn}Every nonempty coarsely homogeneous separable metric space
$X$ of asymptotic dimension zero is macro-uniformly equivalent to one of the next three spaces:
\begin{itemize}
\item $1$ if and only if $X$ is bounded;

\item $2^{<\IN}$ if and only if $X$ is unbounded and has bounded geometry;

\item $\w^{<\IN}$ if and only if $X$  has unbounded geometry.
\end{itemize}
\end{theorem}

This theorem follows from the coarse characterizations of the
Cantor macro-cube and Baire macro-space presented in Theorems
\ref{exc-MU-char} and \ref{Ber-char}. The following coarse
characterization \cite{BZ} of the Cantor macro-cube $2^{<\IN}$ is
an asymptotic analog of the classical Brouwer's characterization
\cite[7.4]{Ke} of the Cantor cube $2^\w$.

\begin{theorem}[Coarse characterization of $2^{<\IN}$]\label{exc-MU-char}
A metric space   $X$ is macro uniformly equivalent to the Cantor macro-space  $2^{<\IN}$ if and only if
\begin{enumerate}
 \item[(1)] $X$ has asymptotic dimension zero;
 \item[(2)] $X$ has bounded geometry;
 \item[(3)] $X$ has no asymptotically isolated balls.
\end{enumerate}
\end{theorem}

Next we present the coarse classification of the Baire macro-space $\w^{<\IN}$. The topological characterization of its topological counterpart $\w^\w$ is a classical result of Aleksandrov and Urysohn (see \cite[7.7]{Ke}): {\em A topological  space  $X$ is homeomorphic to the Baire space  $\w^\w$ if and only if
$X$ is Polish, zero-dimensional and nowhere locally compact.}

\begin{theorem}[Coarse characterization of $\w^{<\IN}$]\label{Ber-char}
A separable metric space  $X$ is macro-uniformly equivalent to the Baire macro-space  $\w^{<\IN}$ if and only if  $X$
has asymptotic dimension zero and has unbounded geometry.
\end{theorem}

We shall prove this theorem in Section 3 using the technique of towers, developed in \cite{BZ}. Now we will look at embeddings of the Baire macro space. First let us recall two classical topological results \cite{Ke}:
\begin{itemize}
\item {\em Each  Polish nowhere locally compact space includes a closed topological copy of the Bare space  $\w^\w$.}
\item {\em Every Polish space is a continuous image of Baire space  $\w^\w$.}
\end{itemize}

There are analogous statements in the coarse category.

\begin{theorem}\label{ob-mk}
Every metric space of unbounded geometry contains a subspace which is macro-uniformly equivalent to the Baire  macro-space  $\w^{<\IN}$.
\end{theorem}

\begin{proof}
Given a metric space $X$ of unbounded geometry, we  have to
construct a macro-uniform embedding $f:\w^{<\IN}\to X$.
Let $\e_1=1$. Taking into account that $X$ has unbounded geometry, by induction construct an increasing unbounded sequence $(\e_i)_{i\in\IN}$ such that $\cov_{6\e_i}^{\e_{i+1}}(X)\ge \w$ for all $i\in\IN$.
For every point $x\in X$ the inequality $\cov_{6\e_i}(B(x,\e_{i+1}))\ge\cov_{6\e_i}^{\e_{i+1}}(X)\ge\w$ implies the existence of a countable subset $S_i\subset B(x,\e_{i+1})$ that contains the point $x$ and is $3\e_i$-separated in the sense that $d(y,z)\ge 3\e_i$ for any distinct points $y,z\in S_i$. Let $f_{i,x}:\w\to S_i$ be any bijective function such that $f_{i,x}(0)=x$. For every $n\in\IN$ let $g_{x,n}:\w^n\to X$ be the function defined by the recursive formula:
$g_{x,1}(i)=f_{1,x}(i)$ and $g_{x,n}(\sigma,i)=g_{f_{n,x}(i),n-1}(\sigma)$ for $\sigma\in\w^{n-1}$. It follows that $f_{x,n}(\sigma,0)=g_{f_{n,x}(0),n-1}(\sigma)=g_{x,n-1}(\sigma)$ for all $\sigma\in\w^{n-1}$. This allows us to define a function $g_{x}:\w^{<\IN}\to X$ letting $g_x|\w^n=g_{x,n}$ for all $n\in\IN$. Here we identify $\w^n$ with the subspace $\{(x_i)_{i\in\IN}\in\w^{<\IN}:\forall i>n\; x_i=0\}$ of $\w^{<\IN}$. One can easily check that the so-defined function $g_x:\w^{<\IN}\to X$ determines a macro-uniform embedding of the Baire macro-space $\w^{<\IN}$ into $X$.
\end{proof}

A metric space $X$ is called {\em macro-connected} if $C_\e(x)=X$ for some $x\in X$ and some $\e<\infty$.
It follows that each unbounded metric space of asymptotic dimension zero is not macro-connected. In particular, the spaces $2^{<\IN}$ and $\w^{<\IN}$ are not macro-connected.

\begin{theorem}\label{image} If a metric space $X$ is not macro-connected, then for each separable metric space $Y$  there is a surjective macro-uniform map $\Phi:X\Ra Y$.
\end{theorem}

\begin{proof} First consider the subspace  $Z=\{n^2:n\in\IN\}$ of the space $\IN$ endowed with the Euclidean metric. Fix any countable dense subset $\{y_n\}_{n=1}^\infty$ in $Y$ and observe that the multi-map $\Phi:Z\Ra Y$, $\Phi:z\mapsto B(y_n,1)$, is macro-uniform and surjective. It remains to construct a surjective macro-uniform map $\psi:X\to Z$.

Fix any points $x_0,x_1\in X$ and let $\e_0=d(x_1,x_0)$. Since $X$ is not macro-connected, there is a sequence of points $(x_i)_{i\in\w}$ of $X$ such that $x_{i+1}\notin C_{\e_i}(x_0)$ where $\e_i=\max\{i,d(x_i,x_0)\}$.

Define a function $\psi:X\to Z$ assigning to each point $x\in X$ the smallest number $n^2\in Z$ such that $x\in C_{\e_n}(x_0)$. It is easy to check that the function $\psi$ is  surjective and macro-uniform.
Then the composition $\Phi\circ\psi:X\Ra Y$ is a required macro-uniform surjective multi-map of $X$ onto $Y$.
\end{proof}

%\section{$\e$-Connected components and uniform multi-maps}

%In this section we recall two lemmas from \cite{BZ}, which  will be used in  proof of Theorem \ref{Ber-char}.
%The following two lemmas can be proved by analogy to Lemmas ?, ? from \cite{BZ}.
%\begin{lemma}\label{l1} Let $\Phi:X\Ra Y$ be a multi-map such that $\Phi^{-1}(Y)=X$.  For any real numbers $\delta\ge 0$ and $\e\ge \w_\Phi(\delta)$, and  every point $x\in X$ the image $\Phi(B_\delta(x))$ lies in the $\e$-connected component $B_\e(y)$ of any point $y\in\Phi(x)$.
%\end{lemma}

%\begin{lemma}\label{l2} Let $\Phi:X\Ra Y$ is a multi-map such that $Y=\Phi(X)$ and $\Phi^{-1}(Y)=X$. For any positive real numbers $\delta<\e$ and $\delta'<\e'$ with $\e'\ge\w_\Phi(\e)$, $\delta\ge\w_{\Phi^{-1}}(\delta')$ we get $\theta_{\delta}^\e(X)\le\theta_{\delta'}^{\e'}(Y)$ and $\Theta_{\delta}^\e(X)\le\Theta_{\delta'}^{\e'}(Y)$.
%\end{lemma}

\section{Towers}

The characterization Theorem~\ref{Ber-char} of the Baire macro-space $\w^{<\IN}$ will
be proved by induction on partially ordered sets called towers. The technique of towers was created in \cite{BZ} for characterization of the Cantor macro-cube $2^{<\IN}$. In this section we recall the necessary information on towers.

\subsection{Partially ordered sets} A {\em partially ordered set} is a set $T$ endowed with a reflexive antisymmetric transitive relation $\le$.

A partially ordered set $T$ is called {\em $\upa$-directed}  if for any two points $x,y\in T$ there is a point $z\in T$ such that $z\ge x$ and $z\ge y$.

A subset $C$ of a partially ordered set $T$ is called  {\em $\upa$-cofinal} if for every $x\in T$ there is $y\in C$ such that  $y\ge x$.

By the {\em lower cone} (resp. {\em upper cone}) of a point $x\in T$  we understand  the set ${\downarrow}x=\{y\in T:y\le x\}$ (resp. ${\uparrow}x=\{y\in T:y\ge x\}$). A subset $A\subset T$ will be called a {\em lower} (resp. {\em upper}) {\em set} if ${\downarrow}a\subset A$ (resp. ${\uparrow}a\subset A$) for all $a\in A$.
For two points $x\le y$ of $T$ the intersection $[x,y]={\upa}x\cap {\da}y$ is called the {\em order interval} with end-points $x,y$.

A partially ordered set $T$ is a {\em tree} if for each point $x\in T$ the lower cone ${\da}x$ is well-ordered (in the sense that each subset $A\subset{\da}x$ has the smallest element).

\subsection{Defining towers}

A partially ordered set $T$ is called a {\em tower} if
$T$ is $\upa$-directed and for every points $x\le y$ in $T$ the order interval $[x,y]\subset T$ is finite and linearly ordered.

This definition implies that for every point $x$ in a tower $T$ the upper set ${\upa}x$ is linearly ordered and is order isomorphic to a subset of $\w$. Since $T$ is $\upa$-directed, for any points $x,y\in T$ the upper sets ${\upa}x$ and ${\upa}y$ have non-empty intersection and this intersection has the smallest element $x\wedge y=\min({\upa}x\cap{\upa}y)$ (because each order interval in $X$ is finite). Thus any two points $x,y$ in a tower have the smallest upper bound $x\wedge y$.

It follows that for each point $x\in T$ of a tower $T$ the lower cone ${\da}x$ endowed with the reverse partial order is a tree of at most countable height.

\subsection{Levels of a tower} Given two points $x,y\in T$ we write $\lev_T(x)\le\lev_T(y)$ if $$|[x,x\wedge y]|\ge|[y,x\wedge y]|.$$ Also we write $\lev_T(x)=\lev_T(y)$ if  $|[x,x\wedge y]|=|[y,x\wedge y]|$.

The relation $$\{(x,y)\in T\times T:\lev_T(x)=\lev_T(y)\}$$ is an equivalence relation  on $T$ dividing the tower $T$ into equivalence classes called the {\em levels} of $T$. The level containing a point $x\in T$ is denoted by $\lev_T(x)$.
Let $$\Lev(T)=\{\lev_T(x):x\in T\}$$ denote the set of levels of $T$ and
$$\lev_T:T\to\Lev(T),\;\lev_T:x\mapsto\lev_T(x),$$
stand for the quotient map called the {\em level map}.

The set $\Lev(T)$ of levels of $T$ endowed with the order $\lev_T(x)\le \lev_T(y)$ is a linearly ordered set, order isomorphic to a subset of integers.
For a level $\lambda\in\Lev(T)$ by $\lambda+1$ (resp. $\lambda-1$) we denote the successor (resp. the predecessor) of $\lambda$ in the level set $\Lev(T)$. If $\lambda$ is a maximal (resp. minimal) level of $T$, then we put $\lambda+1=\emptyset$ (resp. $\lambda-1=\emptyset$).

It is clear that each $\upa$-directed subset $S$ of a tower $T$ is a tower with respect to the partial order inherited from $T$. In this case we say that $S$ is a {\em subtower} of $T$. A typical example of a subtower of $T$ is a {\em level subtower}  $$T^L=\{x\in T:\lev_T(x)\in L\},$$ where $L\subset\Lev(T)$ is an $\upa$-cofinal subset of the level set of the tower $T$.

A tower $T$ will be called {\em $\da$-bounded} (resp. {\em $\upa$-bounded}\/) if the level set $\Lev(T)$ has the smallest (resp. largest) element. Otherwise $T$ is called {\em $\da$-unbounded} (resp. {\em $\upa$-unbounded}\/). In this paper we can consider that all towers are  $\upa$-unbounded and $\da$-bounded.

The level set $\Lev (T)$ of a $\da$-bounded tower can be identified with $\w$, so that zero corresponds to the smallest level of $T$.

\subsection{The boundary of a tower}
By a {\em branch} of a tower $T$ we understand a maximal linearly ordered subset of $T$. The family of all branches of $T$ is denoted by $\partial T$ and is called the {\em boundary} of $T$. The boundary $\partial T$ carries an ultrametric that can be defined as follows.

Given two branches $x,y\in\partial T$ let
$$\rho(x,y)=\begin{cases}0,&\mbox{if $x=y$,}\\
\lev_T(\min x\cap y),&\mbox{if $x\ne y$.}
\end{cases}
$$
It is a standard exercise to check that $\rho$ is a well-defined ultrametric on the boundary $\partial T$ of $T$ turning $\partial T$ into an  ultrametric space.

In the sequel we shall assume that the boundary $\partial T$ of any tower $T$ is endowed with the ultrametric $\rho$.

\subsection{Degrees of points of a tower}

For a point $x\in T$ and a level $\lambda\in\Lev(T)$ let $\suc_\lambda(x)=\lambda\cap{\downarrow}x$ be the set of predecessors of $x$ on the $\lambda$-th level and $\deg_\lambda(x)=|\suc_\lambda(x)|$. For $\lambda=\lev_T(x)-1$,  the set $\suc_{\lambda}(x)$, called the set of parents of $x$, is denoted by $\suc(x)$. The cardinality $|\suc(x)|$ is called the {\em degree} of $x$ and is denoted by $\deg(x)$. Thus $\deg(x)=\deg_{\lev_T(x)-1}(x)$. It follows that $\deg(x)=0$ if and only if $x$ is a minimal element of $T$.

For levels $\lambda,l\in\Lev(T)$ let
$$\deg_\lambda^l(T)=\min\{\deg_\lambda(x):\lev_T(x)=l\}\mbox{ and }\Deg_\lambda^l(T)=\sup\{\deg_\lambda(x):\lev_T(x)=l\}.$$
We shall write $\deg_\lambda(T)$ and $\Deg_\lambda(T)$ instead of $\deg^{\lambda+1}_\lambda(T)$ and $\Deg^{\lambda+1}_{\lambda}(T)$, respectively.

Now let us introduce several notions related to degrees. We define a tower $T$ to be
\begin{itemize}
\item {\em homogeneous} if $\deg_\lambda(T)=\Deg_\lambda(T)$ for any level $\lambda$ of $T$;
\item {\em pruned} if $\deg_\lambda(T)>0$ for every non-minimal level $\lambda$ of $T$.
\end{itemize}

It is easy to check that a tower $T$ is pruned if and only if each branch of $T$ meets each level of $T$.

There is a direct dependence between the degrees of points of the tower $T$ and the capacities of the balls in the ultrametric space $\partial T$. For an arbitrary branch $x\in \partial T$ we can see that $\ent _k (B_n(x))=\deg _k(x\cap \Lev_n(T))$. This implies that $\deg_\lambda^l(T)=\ent_\lambda^l(\partial T)$ and $\Deg_\lambda^l(T)=\Ent_\lambda^l(\partial T)$.

%By a  tower $T_\w$ we understand an ${\upa}$-unbounded homogeneous tower such that  $\deg_\lambda(T)=\w$ for each non-minimal level $\lambda$ of $T$.
%The bare macro space $\w^{<\IN}$  can be identified with the boundary $\partial T_\w$ of a  $\da$-bounded tower $T_\w$.

\subsection{Assigning a tower to a metric space} In the preceding section to each tower $T$ we have assigned the ultrametric space $\partial T$. In this section we describe the converse operation assigning to each metric space $X$ a pruned tower $T_X^L$ whose boundary $\partial T_X^L$ is canonically related to the space $X$.

A closed discrete unbounded subset $L\subset[0,\infty)$ will be called a {\em level set}.
Given a metric space $X$ and a level set $L\subset[0,\infty)$ consider the set $$T^L_X=\{(C_\lambda(x),\lambda):x\in X,\,\lambda\in L\}$$endowed with the partial order $(C_\lambda(x),\lambda)\le (C_l(y),l)$ if $\lambda\le l$ and $C_\lambda(x)\subset C_l(y)$. Here $C_\lambda(x)$ stands for the $\lambda$-connected component of $x$ in $X$.

The tower $T^L_X$ will be called the {\em canonical $L$-tower} of a metric space $X$.
Observe that for each point $x\in X$ the set $C_L(x)=\{(C_\lambda(x),\lambda):\lambda\in L\}$ is a branch of the tower $T^L_X$, so the map
$$C_L:X\to\partial T^L_X,\;\;C_L:x\mapsto C_L(x),$$ called the {\em canonical map}, is well-defined.

The following important fact was proved in \cite[4.6]{BZ}.

\begin{lemma}\label{c5} Let $L\subset[0,\infty)$ be a level set. The canonical map $C_L:X\to\partial T_X^L$ of a metric space $X$ into the boundary of its canonical $L$-tower is a macro-uniform equivalence if and only if $X$ has macro-uniform dimension zero.
\end{lemma}

\subsection{Tower morphisms}
A map $\varphi:S\to T$ is defined to be
\begin{itemize}
\item {\em monotone} if for any $x,y\in S$ the inequality $x<y$ implies $\varphi(x)<\varphi(y)$;
\item {\em level-preserving} if there is an injective map $\varphi_{\Lev}:\Lev(S)\to\Lev(T)$ making the following diagram commutative:
$$
\begin{CD}
S@>{\varphi}>> T\\
@V{\lev_S}VV @ VV{\lev_T}V\\
\Lev(S)@>>{\varphi_{\Lev}}>\Lev(T).
\end{CD}
$$
\end{itemize}

For a monotone level-preserving map $\varphi:S\to T$ the induced map $\varphi_\Lev:\Lev(S)\to\Lev(T)$ is monotone and injective.

A monotone level-preserving map $\varphi:S\to T$ is called
\begin{itemize}
\item {\em a tower isomorphism} if it is bijective;
\item {\em a tower embedding} if it is injective.
\end{itemize}

The following proposition was proved in \cite[5.8]{BZ}.

\begin{proposition}\label{p11} Let $S,T$ be pruned towers and
$f:\Lev(S)\to \Lev(T)$ be a  monotone (and surjective) map. If $\Deg_\lambda^{\lambda+1}(S)\le\deg_{f(\lambda)}^{f(\lambda+1)}(T)$ (and $\deg_\lambda^{\lambda+1}(S)\ge\Deg_{f(\lambda)}^{f(\lambda+1)}(T)$) for each non-maximal level $\lambda\in\Lev(S)$, then there is a tower embedding (a tower isomorphism) $\varphi:S\to T$ such that $\varphi_\lev=f$.
\end{proposition}

Each monotone map $\varphi:S\to T$ between towers induces a multi-map $\partial\varphi:\partial S\Ra\partial T$ assigning to a branch $\beta\subset S$ the set $\partial\varphi(\beta)\subset\partial T$ of all branches of $T$ that contain the linearly ordered subset $\varphi(\beta)$ of $T$. It follows that $\partial\varphi(\beta)\ne\emptyset$ and hence $(\partial\varphi)^{-1}(\partial T)=\partial S$.

\section{Proof of Theorem~\ref{Ber-char}}

To prove the ``only if'' part, assume that a separable metric space $X$ is macro-uniformly equivalent to the Baire macro-space $\w^{<\IN}$ and fix a macro-uniform equivalence $\Phi:X\Ra \w^{<\IN}$. The Baire macro-space $\w^{<\IN}$ is ultrametric and hence has asymptotic dimension zero, see \cite{BDHM}. Since the asymtotic dimension is preserved by macro-uniform equivalences \cite[p.129]{Roe}, the space $X$ also has asymptotic dimension zero. It remains to prove that for every $\delta<\infty$ there is $\e<\infty$ such that $\cov_\delta^\e(X)\ge\omega$. Given $\delta<\infty$, consider the finite number $\delta'=\w_{\Phi}(\delta)$. Since the macro-Baire space $\w^{<\IN}$ has unbounded geometry, there is $\e'<\infty$ such that $\cov_{\delta'}^{\e'}(\w^{<\IN})=\w$. Then for the number $\e=\w_{\Phi^{-1}}(\e')$ we get $\cov_{\delta}^\e(X)\ge\cov_{\delta'}^{\e'}(\w^{<\IN})=\w$.

To prove the ``if'' part, assume that a metric separable space $X$ has asymptotic dimension zero and has unbounded geometry.  Put $\delta_1=1.$ For every natural  $i$, we can find  $\delta_i>\delta_{i-1}+1$ such that
$\ent_{\delta_{i-1}}^{\delta_i}(X)=\w$. Let $L=\{\delta_i\}_{i\in\IN}\subset(0,\infty)$ and consider the canonical
$L$-tower $T_X^L=\{(C_{\lambda}(x),\lambda):x\in X,\;\lambda\in
L\}$ of the metric space $X$. Its level set  $\Lev(T_X^L)$
can identified with set  $L$. By Lemma~\ref{c5}, the
canonical mapping
$$C_L:X\to\partial T_X^L,\;\;C_L:x\mapsto C_L(x)=\{(C_{\lambda}(x),\lambda):\lambda\in L\},$$ is a macro-uniform equivalence.
Since $\cov_{\delta_{i-1}}^{\delta_i}(X)=\Cov_{\delta_{i-1}}^{\delta_i}(X)=\w$,
the tower  $T_X^L$ is homogeneous with $\deg_\lambda(T_X^L)=\Deg_\lambda(T_X^L)=\w$ for each non-minimal $\lambda\in L$.

Let $T_\w$ be the canonical tower of the Baire macro-space $\w^{<\IN}$ with the level set $\IN$. It is clear that $T_\w$ is a homogeneous tower with $\deg_n(T_\w)=\Deg_n(T_\w)=\w$ for each $n\ge 2$.
By Proposition \ref{p11} there is an isomorphism $\varphi:T_X^L\to T_\w$ between the towers $T_X^L$ and $T_\w$. This isomorphism induces a macro-uniform equivalence between the boundaries $\partial T_X^L$ and $\partial T_w=\w^{<\IN}$. Taking into account that $\partial T_X^L$ is macro-uniformly equivalent to $X$, we conclude that $X$ is macro-uniformly equivalent to $\w^{<\IN}=\partial T_\w$.

\section{Proof of Theorem \ref{kl-odn}}
 Let  $X$  be a coarsely homogeneous separable metric space of asymptotic dimension zero.
Since the space $X$ is coarsely homogeneous, there is a function $\varphi:[0,\infty)\to[0,\infty)$ such that for any points $x,y\in X$ there is a macro-uniform equivalence $\Phi:X\Ra X$ such that $y\in\Phi(x)$ and $\max\{\w_{\Phi},\w_{\Phi^{-1}}\}\le\varphi$.

To prove Theorem~\ref{kl-odn}, it is sufficient to check three possible cases.
\smallskip

1. If $X$ is bounded, then the constant map $\Phi:X\to 1=\{0\}$ is a macro-uniform equivalence, so $X$ is coarsely equivalent to the singleton $1$.
\smallskip

2. Now assume that $X$ is unbounded but has bounded geometry. We shall prove that $X$ has no asymptotically isolated balls. Given any $\delta<\infty$ we should find $\e<\infty$ such that $B_\e(x)\ne B_\delta(x)$ for all $x\in X$.

For the number $\delta$ consider the number $\delta'=\varphi(\delta)$. Since the metric space $X$ is unbounded, there are two points $y,z\in X$ on the distance $\e'=d(y,z)>\delta'$. Next, consider the number $\e=\varphi(\e')$. We claim that $B_\e(x)\ne B_\delta(x)$ for all $x\in X$. For this find a macro-uniform equivalence $\Phi:X\Ra X$ such that $y\in\Phi(x)$ and $\w_\Phi\le\varphi$, $\w_{\Phi^{-1}}\le\varphi$. It follows that $\Phi(B_\delta(x))\subset
B_{\w_\Phi(\delta)}(y)\subset B_{\varphi(\delta)}(y)=B_{\delta'}(y)\not\ni z$ and hence $\Phi^{-1}(z)\cap B_\delta(x)=\emptyset$. On the other hand, $\Phi^{-1}(z)\subset \Phi^{-1}(B_{\e'}(y))\subset B_{\w_{\Phi^{-1}}(\e')}(x)\subset B_{\varphi(\e')}(x)\subset B_\e(x)$, which implies that $B_\e(x)\ne B_\delta(x)$.
By Theorem~\ref{exc-MU-char}, the metric space $X$ is macro-uniformly equivalent to the Cantor macro-cube $2^{<\IN}$.
\smallskip

3. Finally, assume that $X$ is not of bounded geometry. Theorem~\ref{Ber-char} will imply that $X$ is macro-uniformly equivalent to the Baire macro-space $\w^{<\IN}$ as soon as we check that $X$ is of unbounded geometry. Assume conversely that $X$ is not of unbounded geometry. This means that there is $\delta<\infty$ such that for every $\e<\infty$ there is a point $x\in X$ with $\cov_\delta(B_\e(x))<\infty$. To derive a contradiction, we shall prove that the metric space $X$ is of bounded geometry. Let $\delta'=\varphi(\delta)$. Given any $\e'<\infty$ we shall prove that $\Cov_{\delta'}^{\e'}(X)<\infty$. Consider the number $\e=\varphi(\e')$ and find a point $x\in X$ such that  $m=\cov_\delta(B_\e(x))<\infty$. We claim that
$\Cov_{\delta'}^{\e'}(X)\le m$. This inequality will follow as soon as we check that  $\cov_{\delta'}(B_{\e'}(y))\le m$ for any point $y\in X$. By the choice of the function $\varphi$, there is a macro-uniform equivalence $\Phi:X\Ra X$ such that $y\in\Phi(x)$ and $\max\{\w_{\Phi},\w_{\Phi^{-1}}\}\le\varphi$. The inequality  $\cov_{\delta}(B_{\e}(x))\le m$ implies the existence of a cover $\U$ of the ball $B_{\e}(x)$ having cardinality $|\U|\le m$ and $\mesh(\U)\le\delta$.

 Then the family $\V=\{\Phi(U):U\in\U\}$ is the cover of the ball $$B_{\e'}(y)\subset\Phi\circ\Phi^{-1}(B_{\e'}(y))\subset \Phi(B_{\w_{\Phi^{-1}}(\e')}(x))\subset \Phi(B_{\varphi(\e')}(x))=\Phi(B_\e(x))$$ and has $\mesh(\V)\le\w_{\Phi}(\delta)\le\varphi(\delta)=\delta'$.
  Since $|\V|\le|\U|\le m$, we conclude that $\cov_{\delta'}(B_{\e'}(y))\le m$. Thus, the space $X$ has bounded geometry and this is a desired contradiction showing that $X$ has unbounded geometry and hence is macro-uniformly equivalent to the Baire macro-space $w^{<\IN}$ according to Theorem~\ref{Ber-char}.

\end{document}